\documentclass{amsproc}

\newtheorem{Theorem}{Theorem}

\newtheorem{theorem}{Theorem}[section]
\newtheorem{lemma}[theorem]{Lemma}
\newtheorem{proposition}[theorem]{Proposition}
\newtheorem{corollary}[theorem]{Corollary}
\theoremstyle{definition}

\theoremstyle{remark}

\numberwithin{equation}{section}

\newcommand{\abs}[1]{\lvert#1\rvert}

\newcommand{\R}{{\mathbb R}}
\newcommand{\Q}{{\mathbb Q}}
\newcommand{\Z}{{\mathbb Z}}

\newcommand{\NN}{{\mathcal N}}

\newcommand{\Int}{{{\rm Int}\,}}

\newcommand{\N}{{\mathbb N}}

\begin{document}

\title[Rotation number]
{Derivatives of rotation number of one parameter families of circle diffeomorphisms}


\author{Shigenori Matsumoto}
\address{Department of Mathematics, College of
Science and Technology, Nihon University, 1-8-14 Kanda, Surugadai,
Chiyoda-ku, Tokyo, 101-8308 Japan
}
\email{matsumo@math.cst.nihon-u.ac.jp
}
\thanks{The author is partially supported by Grant-in-Aid for
Scientific Research (C) No.\ 20540096.}


\date{\today }

\begin{abstract}
We consider the rotation number $\rho(t)$ of a diffeomorphism
$f_t=R_t\circ f$, where $R_t$ is the rotation by $t$ and
$f$ is an orientation preserving $C^\infty$ diffeomorphism of the circle $S^1$.
We shall show that if $\rho(t)$ is irrational
$$\limsup_{t'\to t}(\rho(t')-\rho(t))/(t'-t)\geq 1.$$

\end{abstract}

\maketitle

\begin{center}
\em Dedicated to Professor Kazuo Masuda on his 65-th birthday
\end{center}

\section{Introduction}
Let $f$ be an orientation preserving $C^\infty$ diffeomorphism of the circle $S^1=\R/\Z$,
and consider a one parameter family $f_t$, $t\in J=[-1/2,1/2]$,
 of diffeomorphisms defined by $f_t=R_t\circ f$, where $R_t$ denotes
the rotation by $t$. Fix once and for all a lift $\tilde f:\R\to\R$
of $f$ to the universal covering $\R$ of $S^1$. Then a lift $\tilde f_t$
of
$f_t$ is chosen as $\tilde f_t=T_t\circ \tilde f$, where
$T_t$ is the translation by $t$. The rotation number $\rho(t)\in\R$ of
$\tilde f_t$ is
a continuous and nondecreasing function of $t$. Define a closed set $C$ by
$$
C=J\setminus\Int(\rho^{-1}(\Q)),$$
and assume for simplicity that 
$\rho(-1/2)=0$, $\rho(1/2)$=1 and $C$ is contained in the interior of $J$.

V. I. Arnold \cite{A} showed that $m(C)>0$, where $m$ denotes the
Lebesgue measure.
Denote by $\NN$ the set of non Liouville numbers and define a Borel 
subset $N$ of $C$ by
$$N=\rho^{-1}(\NN).$$
M. R. Herman \cite{H1} showed that $\rho$ is an absolutely continuous
function and that $m(N)>0$, (under a much less restrictive
condition on the one parameter family). A famous theorem of J.-C. Yoccoz says that if $t\in N$,
then $f_t$ is $C^\infty$ conjugate to a rotation.
Thus the result of M. R, Herman says that the set of the value $t$
such that $f_t$ is $C^\infty$ conjugate to a rotation has positive
Lebesgue measure.
On the other hand it is known (\cite{H2} p.170 and \cite{KH} p.412)
that for a generic value $t$ in $C$
the conjugacy of $f_t$ to a rotation is a non absolutely continuous
homeomorphism, provided that $f$ is a real analytic diffeomorphism
and $f'$ is not constantly equal to 1.
Nevertheless it is shown that 
$m(C\setminus N)=0$ (\cite{T}) and furthermore that $\dim_H(C\setminus N)=0$
(\cite{G}), where $\dim_H$ denotes the Hausdorff dimension.
The purpose of this paper is to show a somewhat stronger result
in this direction.

\begin{Theorem} \label{l3}
If $\rho(t)$ is irrational,  then we have
$$
\limsup_{t'\to t}\frac{\rho(t')-\rho(t)}{t'-t}\geq 1.$$
\end{Theorem}

Notice that the above theorem implies by the absolute
continuity of $\rho$ that $\rho^{-1}(B)$ is null
if $B$ is a null Borel set.
As for the case $\rho(t)$ is rational, we have:

\begin{Theorem} \label{t3}
Assume that  $f$ is real analytic and $f'$ is not constantly
equal to 1.
For $t\in C$ such that $\rho(t)\in\Q$,
we have
$$
\limsup_{t'\to t}\frac{\rho(t')-\rho(t)}{t'-t}=\infty.$$
\end{Theorem}

These phenomena can be found in the computer
graphics of the derivative $\rho'$ in \cite{LV}.
The plan of the paper is as follows. In Sect.\ 2,
we prove a weaker version of Theorem \ref{l3} and
apply it to a new proof of the result of \cite{G}.
In Sect.\ 3, we give an ellaboration of the
argument of Sect.\ 2, which yields a proof of
Theorem \ref{l3} for $\rho(t)$ a Liouville number,
while the non Liouville case is treated in
Sect.\ 4.
Finally  Sect.\ 5  is devoted to the proof of Theorem \ref{t3}.
Also we shall remark that it is necessary to consider $\limsup$ instead of 
$\liminf$ in Theorem \ref{l3}.

\section{Weaker version of Theorem \ref{l3}}

The purpose of this section is to show 
the following proposition which is a weaker version of Theorem \ref{l3},
and by applying it to prove that $\dim_H(C\setminus N)=0$ (\cite{G}).

A positive integer $q$
is called a {\em closest return} of an irrational number $\alpha$ if for any
$0<j<q$, we have $\abs{j\alpha}_{S^1}>\abs{q\alpha}_{S^1}$,
where $\abs{x}_{S^1}$ is the distance of $x$ and $0$ in $S^1$.
If $\abs{q\alpha}_{S^1}=\abs{q\alpha-p}$ for some  integer $p$,
the rational number $p/q$ is called a {\em convergent} of $\alpha$.

\begin{proposition} \label{p1}
Suppose $\rho(t)$ is irrational, $p/q$ a convergent of $\rho(t)$,
$t'\in\rho^{-1}(p/q)$ the point nearest to $t$. Then we have
$$
\frac{\rho(t')-\rho(t)}{t'-t}\geq e^{-M},$$
where $M=\Vert (\log f')'\Vert_{C^0}$.
\end{proposition}

To begin with let us prepare the following lemma.

\begin{lemma} \label{l4}
If $\rho(0)$ is irrational, then for any nonnegative integers
$i,j$, we have
$\mu((f^i)'\circ f^j)\geq 1$,
where $\mu$ is the unique $f$-invariant probability measure on
$S^1$.

\end{lemma}

{\sc Proof}. By the downward concavity of log, it suffices to show 
$$
\mu(\log(f^i)'\circ f^j)=0.$$
Since $\mu$ is $f$-invariant, this is equivalent to
$$
\mu(\log(f^i)')=0.$$
Again since $\mu$ is $f$-invariant and
$$\log(f^i)'=\sum_{\nu=0}^{i-1}\log f'\circ f^\nu,$$
this follows from
$$
\mu(\log f')=0.
$$

By the unique ergodicity of $f$, we have a uniform convergence
$$
n^{-1}\log (f^n)'=n^{-1}\sum_{k=0}^{n-1}\log f'\circ f^k
\to \mu(\log f')=a.$$
But if $a>0$ and if $n$ is sufficiently large we have
$$
(f^n)'>\exp{\frac{an}{2}}>1,$$
and if $a<0$, then
$$
(f^n)'<\exp{\frac{an}{2}}<1.$$
In any case these contradict
$$
\int_{S^1}(f^n)'(x)dx=1.$$
\qed

\bigskip

{\sc Proof of Proposition \ref{p1}}. 
Assume that $t$ in Proposition \ref{p1}
is 0. The rotation number of $\tilde f$, $\rho(0)=\alpha$,  
is irrational by the hypothesis.

Given $x,y\in S^1$, denote 
$$
[x,y]=\{z\in S^1\mid x\preceq z\preceq y\},$$
where $\preceq$ is the positive cyclic order of $S^1$, and
by $y-x$ the length of $[x,y]$.

Assume that $p/q$ is a convergent of $\alpha$, 
and, to fix the idea, that  $q\alpha-p<0$. 
Thus we have $\alpha<p/q$, and shall estimate the
 value of $t>0$ such that $\rho(t)=p/q$.

Now since $q$ is a closest return,
the intervals $R_{\alpha}^j[0, -q\alpha+p]$,
$0\leq j\leq q-1$, are mutually disjoint, where
$R_\alpha$ denotes the rotation by $\alpha$. 
The diffeomorphism $f$ is topologically
conjugate to $R_\alpha$ by an orientation preserving homeomorphism
 which maps a given point $x$ to 0.
This implies that if we set $L(x)=[x,f^{-q}(x)]$, then
\begin{equation} \label{e10}
f^jL(x)\ \ \mbox{are mutually disjoint for}\ \ 0\leq j\leq q-1,
\end{equation}
 for any $x\in S^1$.
Let $t$ be the smallest positive value such that $\rho(t)=p/q$. Our aim
is to estimate the value of $(p/q-\alpha)/t$ from below. 

For $0\leq s\leq t$ consider the point $f^q_s(x)$. For $s=0$,
this is just $f^q(x)$ and as $s\to t$ the point $f^q_s(x)$
increases from $f^q(x)$ towards $f^q_t(x)$ on the interval
$$
f^qL(x)=[f^q(x),x].$$
Thus for some $x$, $f^q_t(x)=x$, while for any  $x$
$f^q_t(x)$ lies on $f^qL(x)$, since $t$ is
the smallest value such that $\rho(t)=p/q$.
This shows that the point $f_t^i(x)$ lies on the interval
$f^iL(x)$ for $1\leq i\leq q$. See the figure.
\begin{figure}

\setlength{\unitlength}{1mm}
\begin{picture}(150,75)(10,0)
\put(5,53){$f^q(x)$}
\put(16,53){$f^{q-1}I_1(x)$}
\put(36,53){$f^{q-2}I_2(x)$}
\put(95,53){$I_q(x)$}
\put(110,53){$f^q_t(x)$}
\put(130,53){$x$}
\put(70, 58){$f^qL(x)$}
\put(10,50){\line(1,0){120}}
\put(10,49){\line(0,1){2}}
\put(30,49){\line(0,1){2}}
\put(50,49){\line(0,1){2}}
\put(90,49){\line(0,1){2}}
\put(110,49){\line(0,1){2}}
\put(130,49){\line(0,1){2}}
\put(10,50){\line(1,0){120}}
\put(10,49){\line(0,1){2}}
\put(30,49){\line(0,1){2}}
\put(50,49){\line(0,1){2}}
\put(130,49){\line(0,1){2}}

\put(75,43){$I_{q-1}(x)$}
\put(5, 43){$f^{q-1}(x)$}
\put(10,40){\line(1,0){120}}
\put(10,39){\line(0,1){2}}
\put(30,39){\line(0,1){2}}
\put(50,39){\line(0,1){2}}
\put(90,39){\line(0,1){2}}
\put(70,39){\line(0,1){2}}
\put(130,39){\line(0,1){2}}

\put(5,23){$f^2(x)$}
\put(20,23){$f\circ f_t(x)$}
\put(45,23){$f^2_t(x)$}
\put(34,23){$I_2(x)$}
\put(15,13){$I_1(x)$}
\put(25,13){$f_t(x)$}
\put(5,13){$f(x)$}
\put(10,3){$x$}
\put(130,3){$f^{-q}(x)$}
\put(10,20){\line(1,0){120}}
\put(10,-1){\line(0,1){2}}
\put(30,-1){\line(0,1){2}}
\put(50,-1){\line(0,1){2}}
\put(130,-1){\line(0,1){2}}
\put(10,-1){\line(0,1){2}}
\put(30,-1){\line(0,1){2}}
\put(50,-1){\line(0,1){2}}

\put(130,-1){\line(0,1){2}}

\put(10,10){\line(1,0){120}}
\put(10,19){\line(0,1){2}}
\put(30,19){\line(0,1){2}}
\put(50,19){\line(0,1){2}}
\put(130,19){\line(0,1){2}}
\put(10,19){\line(0,1){2}}
\put(30,19){\line(0,1){2}}
\put(50,19){\line(0,1){2}}
\put(130,19){\line(0,1){2}}

\put(10,0){\line(1,0){120}}
\put(10,9){\line(0,1){2}}
\put(30,9){\line(0,1){2}}
\put(50,9){\line(0,1){2}}
\put(130,9){\line(0,1){2}}
\put(10,9){\line(0,1){2}}
\put(30,9){\line(0,1){2}}
\put(50,9){\line(0,1){2}}
\put(130,19){\line(0,1){2}}

\put(62,5){$f\uparrow$}
\put(62,15){$f\uparrow$}
\put(62,25){$f\uparrow$}
\put(70, 1){$L(x)$}

\put(90,50.5){\line(1,0){20}}
\put(70,40.5){\line(1,0){20}}
\put(30,20.5){\line(1,0){20}}
\put(10,10.5){\line(1,0){20}}
\end{picture}
\bigskip
\end{figure}

For each such $i$, consider the interval
$$I_i(x)=[f\circ f_t^{i-1}(x),f_t^i(x)]\ \ \subset f^iL(x).$$
Since $f^i_t(x)=f\circ f_t^{i-1}(x)+t$, these intervals have length $t$.
Notice that the rightmost point of $I_i(x)$ is mapped by $f$
to the leftmost point of $I_{i+1}(x)$.

The images
$$f^{q-1}I_1(x),f^{q-2}I_2(x)\cdots, fI_{q-1}(x), I_q(x)$$ 
form a sequence of consecutive intervals towards the right.
 Their union 
$$\cup_{i=1}^q f^{q-i}I_i(x)=[f^q(x),f_t^q(x)]$$
is contained in
$f^qL(x)=[f^q(x),x]$ for any $x$.

Let $\tau_i(x)$ be the length of $f^{q-i}I_i(x)$.  
By the Denjoy distortion lemma and (\ref{e10}), we get
$$
\tau_i(x)\geq e^{-M}(f^{q-i})'\circ f^i(x)t,$$
where $M=\Vert (\log f')'\Vert_{C^0}$. Summing them up we get for any $x\in
S^1$,
$$
x-f^q(x)\geq e^{-M}t\sum_{i=1}^q(f^{q-i})'\circ f^i(x).
$$

Now let us evaluate the both hand sides by the invariant measure
$\mu$ of $f$. It is well known that the evaluation of the left term
yields $-1$ times the rotation number  of $\tilde f^q$, i.\ e.\
$$
\mu({\rm Id}-f^q)=p-q\alpha.$$
Therefore Lemma \ref{l4} implies that
$$
(p/q-\alpha)/t\geq e^{-M},$$
as is desired.
\qed

\bigskip
Now let us start the proof of Graczyk's Theorem (\cite{G}). 
First we need the following easy lemma.

\begin{lemma} \label{l100}
Assume $\alpha$ is irrational, $q>1$, and for some $d>3$
$$
\abs{\alpha-p/q}<1/q^{d}.$$
Then $p/q$ is a convergent of $\alpha$.
\qed
\end{lemma}

Let
$$
J_d(p/q)=\rho^{-1}((p/q-1/q^d,p/q+1/q^d))-\rho^{-1}(p/q).$$

\begin{corollary} \label{co}
If $d>3$, $d$ is irrational and $q>1$, then
$$
m(J_d(p/q))\leq 2e^Mq^{-d},$$
where $m$ denotes the Lebesgue measure.
\end{corollary}

{\sc Proof}. Apply Proposition \ref{p1} to the
irrational numbers $\alpha=p/q-1/q^d$ and
$\alpha=p/q+1/q^d$.
\qed

\bigskip
Now the preimage $L$ of the set of Liouville numbers
by $\rho$ can be described as
$$
L=\bigcap_{d>3}\bigcap_{q_0}\bigcup_{q\geq q_0}J_d(p/q),$$
where in the union $p$ runs over the integers $0<p<q$, coprime
to $p$.
Given any $\alpha>0$, if $d>2/\alpha$, we have
$$
\sum_{q\geq q_0}m(J_d(p/q))^\alpha\leq\sum_{q=q_0}^\infty 2^\alpha
e^{\alpha M}q\cdot q^{-\alpha d}\to 0\ \ (q_0\to\infty),$$
which concludes that $\dim_H(L)=0$.

\section{Proof of Theorem \ref{l3} for Liouville $\rho(t)$ }

Let $q_n$ be the $n$-th closest return of the
irrational number $\alpha=\rho(0)$. Then the sequence $\{q_n\alpha\}$
converges to 0 in $S^1$, changing signs alternately. The
closest returns satisfy
$$
q_{n+1}=a_{n+1}q_n+q_{n-1},$$
where $a_{n+1}$ is the $(n+1)$-st denominator of the continued
franction of $\alpha$. 

Here we assume that $\alpha$ is a Liouville number. Thus
the sequence $\{a_{n+1}\}$ is unbounded. It is no loss of generality
to assume that there is a subsequence $n_i$ such that
$$
q_{n_i}\alpha\uparrow 0\ \ {\rm in} \ \ S^1 \ \ {\rm and}\ \ a_{n_i+1}\to\infty.$$
For simplicity we shall write $n_i$ as $n$ in what follows, and have
in mind that 
$$q_n\alpha\prec 0\prec -q_n\alpha\prec q_{n-1}\alpha$$ 
and that $a_{n+1}$ is as large as desired.
All the notations of the previous section are used by replacing
$p/q$ with $p_n/q_n$.

Consider the first return map $S$ of the rotation $R_{\alpha}^{-1}$
on the interval $[q_n\alpha,q_{n-1}\alpha]$. We have
$$
S=\left\{ \begin{array}{ll}
R^{-q_n}_\alpha\ \  {\rm on}\ \ [q_n\alpha,(q_n+q_{n-1})\alpha]
&\mbox{sending it to}\ \ [0, q_{n-1}\alpha]\\
 R_\alpha^{-q_{n-1}}\ \ {\rm on} \ \ 
[(q_n+q_{n-1})\alpha,q_{n-1}\alpha] &\mbox{sending it to}\ \
[q_n\alpha,0].
\end{array}
\right.
$$
Since there are ordering
$$
0\prec (1-a_{n+1})q_n\alpha\prec(q_n+q_{n-1})\alpha\prec -a_{n+1}q_n\alpha
\prec q_{n-1}\alpha,$$
the map $S=R^{-q_n}_\alpha$ sends the interval 
$[0,(1-a_{n+1})q_n\alpha]$ onto $[-q_n\alpha,-a_{n+1}q_n\alpha]$.
In particular $R_\alpha^{-\nu q_n}[0,-q_n\alpha]$, $0\leq \nu<a_{n+1}$,
form a consecutive sequence of intervals contained in $[0, q_{n-1}\alpha]$.
Translating into $f$ via the topological conjugacy, we have for any
$x\in S^1$
\begin{equation} \label{e100}
\bigcup_{\nu=0}^{a_{n+1}-1}f^{-\nu q_n}L(x)\subset K(x),
\end{equation}
where $L(x)=[x, f^{-q_n}(x)]$ as before and $K(x)=[x,f^{q_{n-1}}(x)]$.

As is
well known, easy to show, (\ref{e10}) can be extended to:
\begin{equation} \label{e101}
 f^jK(x)\ \ \mbox{are disjoint for}\ \ 0\leq j\leq q_n-1.
\end{equation}
 So it looks plausible that the total
length $l(x)$ of $\cup_jf^jL(x)$ is very small,
since a large number of its iterates by $f^{-\nu q_n}$ are
mutally disjoint (except the end points), by virtue of
 (\ref{e100}) and (\ref{e101}). 
On the other hand the Denjoy distortion lemma actually guarantees that the
coefficient $e^{-M}$ in Proposition \ref{p1} can be replaced 
by $e^{-Ml}$ where $l$ is the maximum of $l(x)$, and hence 
if $l(x)$ were small enough, we should be able to prove 
Theorem \ref{l3}. 
However we cannot do this for $L(x)$ itself and
instead consider a
subinterval $\hat L(x)$ defined by
$$
\hat L(x)=[x, f^{-q_n}\circ f_t^{q_n}(x)],$$
where as before $t$ is
the smallest value such that $\rho(t)=p_n/q_n$. We are going
to show that the total length of the union
of intervals
$$\hat l(x)=m(\bigcup_{j=1}^{q_n}f^j\hat L(x))$$
is small, where $m$ denotes the Lebesgue measure as before.
 Notice that this is enough for our purpose
of applying the Denjoy distortion lemma, that is, Proposition \ref{p1}
can be improved as
\begin{equation} \label{e55}
\frac{\rho(t)-\rho(0)}{t}\geq e^{-M\hat l},
\end{equation}
where $\hat l=\max \{\hat l(x)\mid x\in S^1\}$.
Now we have
\begin{equation} \label{e102}
\hat L(x)=\bigcup_{i=1}^{q_n}f^{-i}I_i(x),
\end{equation}
where as before
$$I_i(x)=[f\circ f_t^{i-1}(x),f_t^i(x)]\ \ \subset f^iL(x).$$
Put 
$$
\hat l_i(x)=m(\bigcup_{j=1}^{q_n}f^{j-i}I_i(x)).$$
Then we have by (\ref{e102}) and (\ref{e10})
$$
\hat l(x)=\sum_{i=1}^{q_n}\hat l_i(x).
$$
Also (\ref{e100}) and (\ref{e101}) implies that
\begin{equation} \label{e103}
\sum_{\nu=0}^{a_{n+1}-1}\hat l(f^{-\nu q_n}(x))\leq 1.
\end{equation}

Let us compare $\hat l_i(x)$ with
  $\hat l_i(f^{-\nu q_n}(x))$ for $0\leq\nu<a_{n+1}$.
This is possible since the  intervals $I_i(x)$ and
$I_i(f^{-\nu q_n}(x))$ are contained in $f^iK(x)$ and of
length $t$. In fact again by the Denjoy distortion lemma
and (\ref{e101}),
we get
\begin{equation} \label{e105}
m([f^{j-i}I_i(x)])/t\leq e^N m([f^{j-i}I_i(f^{-\nu q_n}(x))])/t
\end{equation}
for any $0<\nu<a_{n+1}$ and $1\leq j\leq q_n$, where
$$
N=\max\{\Vert(\log f')'\Vert_{C^0},\Vert(\log(f^{-1})')'\Vert_{C^0}\}.
$$
Summing up (\ref{e105}) by $j$, we get
$$
\hat l_i(x)\leq e^N\hat l_i(f^{-\nu q_n}x).
$$
Again summing up  by $i$ we obtain
$$
\hat l(x)\leq e^N\hat l(f^{-\nu q_n}(x)).$$
Finally we get by (\ref{e103})
$$a_{n+1}\hat  l(x)\leq e^{N}\sum_{\nu=0}^{a_{n+1}-1}\hat
l(f^{-\nu q_n}(x))
\leq e^N.$$
Now $N$ is a constant depending only on $f$ and $a_{n+1}$
can be chosen arbitrarily large. By virtue of (\ref{e55}),
this completes the proof of
Theorem \ref{l3} for Liouville $\rho(t)$.

\bigskip

\section{Proof of Theorem \ref{l3} for non-Liouville $\rho(t)$ }
 Here we shall show that if
$\rho(t)$ is non Liouville, then $\rho$ is differentiable
at $t$ and  $\rho'(t)\geq 1$.

In the first place we need the following theorem by P. Brunovsk\'y
(\cite{B}).
For the proof see also \cite{H1}.

\begin{theorem} \label{t5}
Let $g_t$ be a $C^1$-path of $C^1$-diffeomorphisms such that
 $g_t$ is 
an irrational rotation for some $t$. Then we have
$$
\frac{d}{dt}{\rm rot}(g_t)=\int_{S^1}\frac{\partial g_t}{\partial t}(x)dx.$$
\end{theorem}

Since $\rho(t)$ is non Liouville, we have $f_t=h\circ R_{\rho(t)}\circ
h^{-1}$
for some $C^\infty$ diffeomorphism $h$ (\cite{Y}). 
Applying the Brunovsk\'y theorem
to the family $h^{-1}\circ f_t\circ h$, we get
\begin{eqnarray*}
&\rho'(t)=\int_{S^1}\frac{\partial (h^{-1}\circ f_t\circ h)}
{\partial t}(x)dx
=
\int_{S^1}(h^{-1})'\circ f_t\circ h(x)\frac{\partial f_t}
{\partial t}\circ h(x)dx\\
&=
\int_{S^1}(h^{-1})'\circ f_t\circ h(x)dx.
\end{eqnarray*}

Since $f_t\circ h=h\circ R_{\rho(t)}$ and the Lebesgue measure is
invariant by the rotation we have
$$
\rho'(t)=\int_{S^1}(h^{-1})'\circ h(x)dx=\int_{S^1}h'(x)^{-1}dx.$$
Now the Schwarz inequality 
concludes the proof of Theorem \ref{l3} for non-Liouville $\rho(t)$.

\section{Proof of Theorem \ref{t3}.}
By the assumption of Theorem \ref{t3}, for any $p/q\in\Q$,
the set $\rho^{-1}(p/q)$ is a nondegenerate interval and $C$
is a Cantor set. It is no loss of generality to assume that
$0\in J$ is the supremum of $\rho^{-1}(p/q)$ and to show
$$
\lim_{t\downarrow 0}\frac{\rho(t)-(p/q)}{t}=\infty.$$

The real analyticity of $f$ implies that
the periodic points of $f_0=f$ are finite in number, say $x_\nu$,
$1\leq \nu\leq ql$. Now since $0$ is the
supremum
of $\rho^{-1}(p/q)$,
the graph of $\tilde f^q-p$ is above the diagonal and tangent to it
just at the points $\pi^{-1}(x_\nu)$. For $\varepsilon>0$, let
$I_\nu=[a_\nu,b_\nu]$ be the $\varepsilon$-neighbourhood of $x_\nu$.
Choosing $\varepsilon$ small enough, one can assume that the intervals
$I_\nu$ are disjoint. Put $J_\nu=[b_\nu,a_{\nu+1}]$.
Choose $N>0$ big enough so that 
$$
f^{qj(\nu)}(b_\nu)\not\in J_\nu, \ \ 1\leq \exists j(\nu)\leq N+1. 
$$
This means that any orbit by $f^q$ stays consecutively in $J_\nu$
for at most $N$ times.
Then since $\tilde f_t^q>\tilde f^q$ for $t>0$ and the speed
for $f_t^q$ is bigger than that for $f^q$,
any orbit by $f_t^q$ stays consecutively in $J_\nu$
for at most $N$ times.

On the other hand straightforward computation shows that 
$\partial f_t^q/\partial 
t\geq 1$, and therefore $\tilde f_t^q-p\geq t$.
This shows that any orbit by $f_t^q$ stays consecutively in $I_\nu$ for
at most $2\varepsilon/t$ times, $2\varepsilon$ being the length of $I_\nu$.

Let us estimates the times of iterations
of $f_t^q$ needed for some
point to go around $S^1$ once. 
The above observation shows that the times needed 
for a round trip does not exceed $(N+2\varepsilon/t)ql$.
Translated into  the language of rotation number we have
$$
q\rho(t)-p\geq ((N+2\varepsilon/t)ql)^{-1}.$$
Therefor if $t<\varepsilon/N$, we have
$$\frac{\rho(t)-(p/q)}{t}\geq (3\varepsilon q^2l)^{-1},
$$
completing the proof of Theorem \ref{t3}.

\bigskip
Finally let us remark that taking $\limsup$ instead of $\liminf$
is necessary for Theorem \ref{l3}. 
\begin{proposition}
Assume that  $f$ is real analytic and $f'$ is not constantly equal to 1.
There is a residual subset $R$ in $C$ such that for any $t\in R$
$$
\liminf_{t'\to t}\frac{\rho(t')-\rho(t)}{t'-t}=0, \ \
{\rm and}\ \ \limsup_{t'\to t}\frac{\rho(t')-\rho(t)}{t'-t}=\infty.$$
\end{proposition}

{\sc Proof}. Let $\rho(J)\cap\Q=\{\alpha_n\mid n\in\N\}$ and
set $[a_n,b_n]=\rho^{-1}(\alpha_n)$. Then there is $c_n>b_n$ very near $b_n$
such that if $t\in(b_n,c_n)$
$$
\frac{\rho(t)-\alpha_n}{t-a_n}<\frac{1}{n}\ \ {\rm and}
\ \ \frac{\rho(t)-\alpha_n}{t-b_n}>n.
$$

We used Theorem \ref{t3} for the second inequality.
Now the set
$$
R=\cap_p\cup_{n>p}(b_n,c_n)\cap C$$
is residual in $C$ and satisfies the condition of the theorem.
\qed


\begin{thebibliography}{99}
\bibitem[A]{A} V. I. Arnold, {\em Small denominator I:
on the mapping of a circle into itself}, Isvestija Akad.\ Nauk.\
serie Math., {\bf 25}, 1(1961), 21-86,  Translations A. M. S.
2nd Series {\bf 46}, 213-284.

\bibitem[B]{B} P. Brunovsk\'y, {\em Generic properties of
the rotation number of one-parameter diffeomorphisms of the circle},
Czech.\ Math.\ J. {\bf 24}(1974), 74-90.

\bibitem[G]{G} J. Graczyk, {\em Linearizable circle diffeomorphisms in
one-parameter families,} Bol.\ Soc.\ Bras.\ Mat.\ {\bf 24}(1993),
201-210. 


\bibitem[H1]{H1} M. R. Herman, {\em Mesure de Lebesgue et nombre de rotation},
Springer Lecture Notes in Math., {\bf 597}(1977), 271-293.

\bibitem[H2]{H2} M. R. Herman, {\em Sur la conjugaison diff'erentiable
des diff\'eomorphismes du cercle a des rotations,} Publ.\ I. H. E. S. 
{\bf 49}(1979), 5-233.

\bibitem[KH]{KH} A. Katok and B. Hasselblatt, {\em Introduction
to the Modern Theory of Dynamical Systems,} Cambridge University
Press, 1995.

\bibitem[LV]{LV} A. Luque and J. Villanueva, {\em Computation of
derivatives of the rotation number for parametric families of
circle diffeomorphisms,} Physica D {\bf 237}(2008) 2599-2615.


\bibitem[T]{T} M. Tsujii, {\em Rotation number and one-parameter 
families of circle diffeomorphisms,} Erg.\ Th.\ Dyn.\ Sys.\
{\bf 12}(1992), 359-363.

\bibitem[Y]{Y} J.-C. Yoccoz, {\em Conjugaison diff\'erentiable des 
diff\'eimorphismes du cercle dont le nombre de rotation v\'erifie une
conditon diophantinne,}
Ann.\ Sci.\ Ecole Norm.\ Sup.\ {\bf 17}(1984), no.\ 3, 333-359.


\end{thebibliography}
\end{document}